\documentclass[a4paper, 11pt]{amsart}
\usepackage[spanish]{babel}
\input xy   
\xyoption{all}
\swapnumbers

\begin{document}

\title{Categor\'{i}as. Los 30 primeros a\~{n}os}
\author{Eduardo J. Dubuc} 

\maketitle

El lector debe tener en cuenta que esta recopilaci\'on de los 30 primeros a\~{n}os de vida de la teor\'{\i}a de categor\'{\i}as responde a lo que puedo recordar, y no es el resultado de una investigaci\'{o}n hist\'orica sobre el desarrollo de la misma. El relato tiene como hilo conductor las veinti\'un referencias primarias que aparecen por orden cronol\'ogico al final del art\'{\i}culo, las secundarias ir\'an en notas a pie de p\'agina.  

Asimismo, he querido que estas notas sean un fiel reflejo de la conferencia impartida en octubre de 2013 en la ENHEM 4\footnote{Cuarta Escuela Nacional de Historia y Educaci\'on Matem\' 
atica, Cali (Colombia).}, de contenido acotado por el tiempo, y no incluyen ni m\'as ni menos que lo all\'{\i} dicho. Quiero agradecer a los organizadores de la ENHEM 4, y en particular a Luis Recalde, el haberme invitado a la Escuela, sin ello este trabajo nunca hubiese visto la luz. 

\section*{La definici\'on de categor\'{i}a y los primeros pasos}

{\bf[1941]} 
Hacia 1941, en la teor\'{i}a de la homotop\'{i}a se empez\'o a representar funciones como flechas y ecuaciones como diagramas conmutativos. 
\begin{equation} \label{diagrama}
\xymatrix
        {
         h=gf
        }
\hspace{20ex}
\xymatrix@R=10pt@C=7pt
         {
          X \ar[rr]^f  \ar[rd]_h   &  &  Y \ar[ld]^g
           \\
            & Z  
         } 
\end{equation}
La validez de la ecuaci\'on dice que $h$ es la \emph{composici\'on} de $f$ con $g$, y se expresa gr\'aficamente diciendo que el correspondiente tri\'angulo es \emph{commutativo}.

\vspace{1ex}

Esta notaci\'on condujo a un concepto: se tienen \emph{objetos}, se tienen \emph{flechas}, y se tienen diagramas que expresan la composici\'on de flechas. Las propiedades b\'asicas son las ecuaciones de \emph{identidad}, $f \, id_X = f = id_Y f$, y la \emph{asociatividad} de la composici\'on, $h(gf) = (hg)f$. Es decir, se tiene una \emph{categor\'{i}a}.

\vspace{1ex} 

{\bf[1942-1945]} La definici\'on de categor\'{i}a apareci\'o en el trabajo de Eilenberg-Mac Lane\footnote{Samuel Eilenberg (1913-1998), Sanders Mac Lane (1909-2005).} \cite[1945]{EM2} para dilucidar la estructura abstracta subyacente en los desarrollos de \cite[1942]{EM1}. En una categor\'{\i}a hay una correspondencia biyectiva entre los morfismos identidad ``$id_X$'' y los objetos ``$X$'', de forma tal que estos pueden ser omitidos en la definici\'on de categor\'{\i}a. Sin embargo, Eilenberg y Mac Lane escriben:
\begin{quote} \label{objetos}
\emph{Es claro por tanto que los objetos juegan un papel secundario, y podr\'{\i}an ser completamente omitidos en la definici\'on de categor\'{\i}a. Sin embargo, la manipulaci\'on de las aplicaciones ser\'{\i}a algo menos conveniente si as\'{\i} se hiciera.}\footnote{``It is thus clear that the objects play a secondary role, and could be entirely omitted from the definition of a category. However, the manipulation of the applications would be slightly less convenient were this done''.  \cite[p\'ag. 238]{EM2}} 
\end{quote} 
 As\'{\i}, desde el comienzo, una clara distinci\'on entre objetos y flechas, acorde con la notaci\'on de diagramas, hace a la esencia misma de la noci\'on de categor\'{\i}a, aunque desde un punto de vista puramente formal los objetos podr\'{\i}an haber sido completamente dejados de lado.

En \cite{EM2} se establece expl\'{\i}citamente el principio: 
\emph{Dado cualquier objeto matem\'atico debe considerarse tambi\'en una noci\'on de morfismo.} Se llega as\'{i} a la sucesi\'on de nociones: 
\begin{center}
     \emph{Categor\'{i}as $\;\leadsto\;$ Funtores  $\;\leadsto\;$ Transformaciones naturales}.
\end{center}

En los primeros ejemplos este lenguage se utiliz\'o para clarificar y definir en forma compacta un concepto de gran complejidad, como es la comparaci\'on \mbox{(\emph{transformaci\'on natural $\alpha$})} entre dos teor\'{\i}as de homolog\'{\i}a (\emph{funtores $F$ y $G$}) que asocian a cada espacio $X$ 
(\emph{objeto de la categor\'{i}a $\mathcal{G}$}) un \'algebra \mbox{(\emph{$FX$, $GX$, objetos de la categor\'{i}a 
$\mathcal{A}$}).}
Todo esto se resume con una representaci\'on conjunta en el sencillo diagrama
\begin{equation} \label{Cat}
\mathcal{G}
\xymatrix@C=9ex@R=2.4ex
    {
     {} \ar@<1.6ex>[r]^{F}
        \ar@{}@<-1.4ex>[r]^{\!\!\alpha\,\Downarrow}
        \ar@<-1.3ex>[r]_{G}
    & {}
    }
\mathcal{A}\,
\hspace{8ex}
FX \overset {\alpha_X} {\longrightarrow} GX
\end{equation}

En estas obras la teor\'{i}a de categor\'{i}as es desarrollada y utilizada esencialmente como una notaci\'on y un lenguaje, y tambi\'en para unificar construcciones. Un ejemplo considerado expl\'{\i}citamente en \cite{EM2} es la noci\'on de \emph{l\'{i}mite inverso}, que unifica la construcci\'on de los  n\'umeros $p$-\'adicos  con la construcci\'on de la homolog\'{i}a de C\u{e}ch:

\vspace{1ex}

Dado el conjunto ordenado de los n\'umeros 
$\{\; \cdots\; \geq n\geq n-1\geq \;\cdots \; \geq 0\}$, se tiene el sistema inverso determinado por el grupo de enteros m\'odulo $p^n$, cuyo l\'{\i}mite inverso es el grupo de enteros $p$-\'adicos:
$$
\mathbb{Z}_p \;\;\; \to \;\;\; 
\left.\rule{0pt}{2.5ex}\right \{ 
\;\; \cdots\;\;\to\frac{\mathbb{Z}}{p^n\mathbb{Z}}\to\frac{\mathbb{Z}}{p^{n-1}\mathbb{Z}}\to\;\; \cdots\;\;\to\{0\}
\left.\rule{0pt}{2.5ex}\right \}
$$

Un espacio topol\'ogico $X$ determina el conjunto de cubrimientos abiertos ordenado por refinamiento 
$\{ \;\; \cdots \;\; \succ \mathcal{V}\succ\mathcal{U}\succ \;\; \cdots\;\; \succ\{X\} \}$. Dado un grupo $G$ y un entero $q$, se tiene el sistema inverso determinado por el $q$-\'{e}simo grupo de homolog\'{i}a del nervio del cubrimiento con coeficientes en $G$, cuyo l\'{\i}mite inverso es el $q$-\'{e}simo grupo de homolog\'{\i}a de 
C\u{e}ch de $X$:
$$
\text{C\u{e}ch}(X) \;\;\; \to \;\;\; 
\left.\rule{0pt}{2.5ex}\right \{ 
\;\; \cdots\;\;H^q(N(\mathcal{V}),G)\to H^q(N(\mathcal{U}),G)\to\;\; \cdots\;\;\to\{0\}
\left.\rule{0pt}{2.5ex}\right \}
$$

En ese primer momento, estos l\'{\i}mites se construyen con los elementos de los conjuntos subyacentes. No se considera la \emph{propiedad universal} del l\'{i}mite. 

\section*{La noci\'on de propiedad universal}

{\bf [1948]} Parece que la noci\'on de propiedad universal surge por primera vez en 1948, en trabajos sobre topolog\'{\i}a de Pierre Samuel, que era miembro de Bourbaki. Las propiedades universales fueron incorporadas por \mbox{Bourbaki} en su libro sobre la teor\'{\i}a de conjuntos y el concepto de estructura \mbox{\cite[1939-1957]{B}.} All\'{\i} considera, por ejemplo, el grupo libre $L(S)$ sobre un conjunto $S$ y el producto tensorial $G \otimes H$ de dos grupos abelianos, expresadas aqu\'{\i} abajo por sus diagramas:
$$
\xymatrix@R=5ex
         {                                       
           S  \ar@/_/[drr]_{\forall\;\text{funci\'on} \;} \ar[rr]^\eta   &&   L(S)   \ar@{-->}[d]^{\exists! \;\text{morfismo}}    
          \\
                                       &&   G   
         }
\hspace{1cm}
\xymatrix@R=5ex
         {                                       
           G\times H  \ar@/_/[drr]_{\forall\;\text{bilineal} \;\;\;\;} \ar[rr]^{\text{bilineal}}   &&   G\otimes H   \ar@{-->}[d]^{\exists ! \;\text{morfismo}}   
          \\
                                       &&   B   
         }
$$ 

All\'{\i} tambi\'en se consideran las construcciones de l\'{\i}mite inverso (o \emph{proyectivo}) y de \emph{l\'{\i}mite directo} (o \emph{inductivo}). Estas son construcciones totalmente distintas, pero ambas caracterizadas por propiedades universales, y, significativamente, esta caracterizaci\'on es el primer enunciado demostrado por Bourbaki inmediatamente despu\'es de las respectivas definiciones conjuntistas. Luego figuran una serie de proposiciones todas demostradas a partir de la propiedad universal, y no de las construcciones expl\'{\i}citas de la definici\'on. Sin embargo, Bourbaki se ve obligado a escribir en detalle dos veces cada demostraci\'on, una para cada concepto. Ello es debido a que, si bien se trata de propiedades universales duales una de la otra (por lo que basta una \'unica demostraci\'on para validar ambas), para el propio enunciado de este hecho es necesario contar con la noci\'on de categor\'{\i}a, de propiedad universal definida con ese lenguage, y del \mbox{principio de dualidad de Mac Lane considerado m\'as abajo.} 

Cuenta Armand Borel (1923-2003)\footnote{Armand Borel, ``Twenty-five years with Nicolas Bourbaki (1949-1973)'',  \textit{Notices of the AMS} {\bf 45}(3) 1998, 373--380.} que Alexander Grothendieck hab\'{\i}a presentado a Bourbaki un proyecto (detallado) para escribir los fundamentos en el lenguaje categ\'orico, pero fue descartado luego de largas discusiones, porque finalmente prim\'o la pol\'{i}tica que afirmaba que generalizaciones que no agregan aplicaciones deben descartarse, versus la pol\'{i}tica de Grothendieck que sosten\'{\i}a que los fundamentos deber\'{\i}an estar hechos  \emph{no solo para la matem\'atica existente sino tambi\'en para la posible matem\'atica futura.}\footnote{En palabras de Borel: ``his plan aimed at supplying foundations not just for existing mathematics, as had been the case for the ``Elements'', but also for future developments to the extent they could be foreseen''. [\emph{op. cit}. p\'ag. 378]} 

Es posible pensar que tambi\'en haya influido en la decisi\'on de Bourbaki el hecho que aceptar el proyecto de Grothendieck implicaba una reescritura completa de todo un libro que ya estaba escrito en su forma final. Tambi\'en se aduce (en mi opini\'on equivocadamente)  que utilizar las categor\'{\i}as implicaba dificultades debido a que la teor\'{\i}a de conjuntos no resultaba suficiente para definir la noci\'on de categor\'{i}a en forma adecuada. Sin embargo, estas dificultades son puramente de naturaleza t\'ecnica, y ya hab\'{\i}an sido satisfactoriamente resueltas por Grothendieck con su teor\'{\i}a de universos (cardinales fuertemente inaccesibles), m\'as tarde publicada en un ap\'endice de \cite{SGA4} \mbox{firmado por Bourbaki.}

\section*{La noci\'on categ\'orica de propiedad universal.}

{\bf [1948-1950]} Pas\'o una d\'ecada sin grandes novedades en la propia teor\'{i}a de categor\'{i}as. Los inventores del concepto trabajaron con enunciados categ\'oricos pero en categor\'{i}as particulares. Mac Lane lo hizo en art\'{i}culos sobre categor\'{i}as de grupos 
\mbox{\cite[1948, 1950]{ML1,ML2}.} All\'{\i} utiliza la riqueza del lenguaje categ\'orico, que permite considerar \emph{enunciados duales} (el enunciado dual se obtiene invirtiendo la direccion de las flechas, como est\'a ilustrado aqu\'{\i} abajo) y \emph{demostraciones duales}, lo que resulta posible por ser los axiomas de categor\'{i}a \emph{autoduales}. Expl\'{\i}citamente establece el principio de dualidad: \emph{Si un enunciado sobre una categor\'{\i}a se sigue de los axiomas de categor\'{\i}a, el enunciado dual tambi\'en se sigue}\footnote{``DUALITY PRINCIPLE. If any statement about a category is deducible from the axioms for a category, the dual statement is likewise deducible''. \cite[p\'ag. 498]{ML2}}. Tambi\'en introduce la definici\'on de objetos por medio de propiedades universales en lugar de las cl\'asicas construcciones conjuntistas. Esto permite autom\'aticamente la consideraci\'on de las definiciones duales. 

\medskip

Como ejemplos de definiciones universales y sus duales tenemos las siguientes, expresadas por sus diagramas (entre par\'entesis indicamos la correspondiente construcci\'on en la categor\'{\i}a de los conjuntos). En la primera fila, el \emph{producto cartesiano} ``$\,\times\,$'' (conjunto de pares de elementos) y el \emph{objeto final} ``$\,1\,$'' (conjunto de un solo elemento $\{*\}$), y debajo, en la segunda, sus respectivos duales, el \emph{coproducto} ``$\,\amalg\,$'' (uni\'on disjunta) y el \emph{objeto inicial} ``$\,0\,$'' (conjunto sin elementos $\emptyset$):
$$
\xymatrix@R=2ex
         {
                                       &&    A  
           \\
           X \ar@/^/[urr] \ar@/_/[drr] \ar@{-->}[rr]^{\exists !}   &&   A\times B \ar[u]  \ar[d]    
          \\
                                       &&   B 
         }
\hspace{1cm}  \xymatrix@R=2ex
         {
                                       &&   
           \\
           X  \ar@{-->}[rr]^{\exists !}   &&   1  
          \\
                                       &&  
         }
$$

$$
\xymatrix@R=2ex
         {
                                       &&    A   \ar[d]   \ar@/_/[dll]
           \\
           X     &&   A \amalg B    \ar@{-->}[ll]_{\exists !} 
          \\
                                       &&   B \ar[u] \ar@/^/[ull]
         }
\hspace{1cm}  \xymatrix@R=2ex
         {
                                       &&   
           \\
           X     &&   0  \ar@{-->}_{\exists !}[ll]
          \\
                                       &&  
         }
$$

Este ejemplo, pero en la categor\'{\i}a de los grupos, es considerado expl\'{\i}citamente en \cite{ML2}. En ese caso, ``$\,\times\,$'' sigue siendo el grupo de pares de elementos, y ``$\,1\,$'', el grupo de un solo elemento $\{e\}$\footnote{Se usa la letra $e$ para denotar la unidad de un grupo.}. Pero ``$\,\amalg\,$'' no es la uni\'on disjunta (que no tiene estructura de grupo) sino que es el producto libre de grupos.
``$\,0\,$'' es ahora tambi\'en $\{e\}$, as\'{\i}, $1 = 0$ en la categor\'{\i}a de grupos.

Otro ejemplo considerado lo proporciona el grupo cociente asociado a un subgrupo $K\subseteq G$, que da lugar a un homomorfismo $\rho: G \longrightarrow G/K$ que verifica $\rho(K)=e$, donde el grupo $G/K$ es un conjunto de clases de equivalencia. Ahora, la noci\'on se enuncia como la propiedad universal del homomorfismo $\rho$ entre todos los que tienen la propiedad $\varphi(K)=e$: 
$$
\xymatrix@R=5ex
         {                                       
           G  \ar@/_/[drr]_{\forall \, \varphi} \ar[rr]^\rho   &&   G/K   \ar@{-->}[d]^{\exists !}  & \rho(K)=e   
          \\
                                       &&   B   &  \varphi(K)=e
         }
$$ 

Definidos de esta manera, los objetos son caracterizados salvo isomorfismo (\'unico). As\'{\i}, cualquier suryecci\'on de grupos $\varphi: G \longrightarrow H$ puede ser considerada el cociente $H=G/K$ por el subgrupo $K = \text{N\'ucleo}(\varphi)=\varphi^{-1}(e)$. Esta manera de pensar ``categ\'orica'' es tan radicalmente diferente de la cl\'asica manera ``conjuntista'', que el mismo Mac Lane en ese entonces todav\'{\i}a se le acerca con precauci\'on, t\'{\i}midamente. Por ejemplo, en el mismo trabajo, en su tratamiento de la dualidad observa que contrariamente a la relaci\'on de subgrupo, la relaci\'on dual, de cociente, no es transitiva pues un cociente de un cociente no es un cociente (pensando cociente en el sentido conjuntista de clases de equivalencia), y hace toda una complicada elaboraci\'on con congruencias para arreglar el asunto. Sin embargo, de haber definido cociente por propiedad universal, la relaci\'on resulta autom\'aticamente transitiva dado que cualquier suryecci\'on entre grupos es un cociente.

\section*{Dos libros fundamentales.}

{\bf[1950, 1953/56]}  En estos a\~{n}os Eilenberg, en colaboraciones con N. Steenrod (1910-1971) y con H. Cartan (1904-2008), publica dos libros fundamentales, \cite[1950]{ES} y \cite[1956]{CE}\footnote{Seg\'un Cartan, el segundo libro fue entregado a Princeton University Press en 1953. V\'ease Hyman Bass, Henri Cartan, Peter Freyd, Alex Heller, and Saunders Mac Lane, ``Samuel Eilemberg 1913-1998'', \emph{Bibliographical Memoirs} Vol. 79, The National Academic Press 2000.}, que hacen uso extensivo de las categor\'{\i}as, funtores y transformaciones naturales como una eficiente herramienta para demostraciones en la categor\'{\i}a de m\'odulos. Estos libros hicieron \'epoca, introdujeron el t\'ermino \emph{\'algebra homol\'ogica}, y dejaron la matem\'atica madura para el desarrollo de la teor\'{\i}a de categor\'{\i}as en s\'{\i} misma. Sin embargo, todav\'{\i}a no hacen uso del concepto de propiedad universal.

\section*{El lema de Yoneda}

Las categor\'{i}as son la estructura apropiada para tratar la noci\'on de propiedad universal en su generalidad correcta, y a su vez esta noci\'on hace que la teor\'{\i}a de categor\'{\i}as sea algo m\'as que un exitoso lenguage (para establecer enunciados) y una poderosa herramienta (para establecer demostraciones) como hasta el momento hab\'{\i}a sido utilizada exitosamente en los libros \cite{ES} y \cite{CE} antes citados. 

Las definiciones por \emph{propiedad universal} implican una revoluci\'on en el pensamiento matem\'atico. Su \'{\i}ntima relaci\'on con la noci\'on de \emph{funtor representable} dada por el \emph{lema de Yoneda} es lo que hace de las categor\'{i}as una teor\'{i}a con construcciones y teoremas profundos. 

\vspace{1ex}

{\bf[1954]}  Mac Lane y N. Yoneda (1930-1996) coincidieron en Par\'{\i}s al principio de los a\~{n}os 50, donde Mac Lane aprendi\'o de boca del propio Yoneda sus trabajos en \'algebra homol\'ogica. Mac Lane identific\'o, oculto como un caso particular en el c\'omputo de los funtores \emph{Ext} en categor\'{\i}as de m\'odulos, un resultado al que le otorg\'o una importancia que iba mucho m\'as all\'a que el c\'alculo para el que lo utiliz\'o Yoneda. As\'{\i} lo populariz\'o en varias charlas (anteriores a la publicaci\'on de Yoneda \cite[1954]{Y}), y lo bautiz\'o con el nombre ``Lemma de Yoneda'', enunci\'andolo ya en categor\'{\i}as abstractas\footnote{http://math.stackexchange.com/questions/53656/what-is-the-origin-of-the-expression-yoneda-lemma}. 

\vspace{1ex}

{\bf Propiedades universales y funtores representables.}
Un funtor $F: \mathcal{X} \to \mathcal{Y}$ entre categor\'{\i}as es una funci\'on que asigna objetos a objetos y flechas a flechas, respetando la estructura, lo que significa que env\'{\i}a tri\'angulos commutativos en tri\'angulos commutativos. El funtor se dice \emph{contravariante} si invierte la direcci\'on de las flechas, lo que  se indica con ``op'', $F: \mathcal{X}^{op} \to \mathcal{Y}$. Dada una categor\'{i}a $\mathcal{X}$ y en ella un objeto $C$, se tiene el funtor contravariante a valores en los conjuntos, 
\mbox{$hom (-,C):\mathcal{X}^{op}\to Set$,} definido en los objetos por el conjunto de flechas $hom (X,C) = \{f | f: X \to C \}$. Un funtor contravariante $F: \mathcal{X}^{op}\to Set$ es \mbox{\emph{representable}} si es isomorfo a un \mbox{funtor $hom(-,C)$.} 

\medskip

Yoneda \cite[1954]{Y} estableci\'o el resultado conocido cono \emph{Lema de Yoneda}. Dada una transformacion natural 
$\theta:hom(-,C)\to F$, el lema afirma que la asignaci\'on del elemento $\xi=\theta_C(id_C) \in F(C)$ determina una biyecci\'on entre el conjunto de transformaciones naturales $\{hom(-,\, C) \to F \}$ y el conjunto $F(C)$, $Nat(hom(-, C), F) \cong F(C)$. Esta biyecci\'on usualmente se describe con el siguiente diagrama (donde indicamos tambi\'en la aplicaci\'on inversa):
$$
\frac{\theta_X:hom(X,C)\to F(X), \hspace{2ex}\theta_X(f)=F(f)(\xi)}{\xi\in F(C), \hspace{2ex}\xi=\theta_C(id_C)}.
$$
Que $F$ sea un funtor representable significa entonces que se tiene un objeto $C$ y un elemento $\xi \in F(C)$ tal que la correspondiente transformaci\'on  $\theta_X$ es biyectiva para todo $X$. Es decir:
\begin{itemize}
\item[(i)] $\xi\in F(C)$.
\item[(ii)] $\forall \, x\in F(X)$, $\exists! f:X\to C 
\, \mid \, F(f)(\xi)=x$.
\end{itemize}

El objeto $C$ (junto con $\xi$) es la soluci\'on a una \emph{propiedad universal} en la categor\'{\i}a $\mathcal{X}$. Toda propiedad universal corresponde a un funtor $\text{Datos}:\mathcal{X}^{op} \to Set.$ La existencia de soluci\'on significa que se tiene un \emph{dato universal}, es decir, que el funtor ``Datos'' es representable. Notar que nos estamos mordiendo la cola:

 \emph{La definici\'on abstracta de propiedad universal en una categor\'{\i}a 
$\mathcal{X}$ es, por definici\'on, la representabilidad de un funtor definido en $\mathcal{X}$ a valores en los conjuntos.}

\vspace{1ex}

{\bf Moduli spaces.} Son de este tipo los problemas llamados ``de Moduli'' (Moduli spaces), versi\'on actual del cl\'asico problema que consiste en darle una estructura geom\'etrica al conjunto de objetos geom\'etricos que se quieren clasificar y estudiar c\'omo var\'{i}an (modulan) continuamente (familias continuas parametrizadas por un espacio).
\begin{figure}
$$
\xymatrix@R=1.5ex@C=0.2ex
       { 
         {}&{}&{}&{}&{}
        &{\cdot}\ar@{-}[dddd]
        &{}&{}&{}&{}&{}&{}&{}&{}&{}&{}
        &{}&{}&{}&{}&{}&{}&{}&{}
        &{\cdot}\ar@{-}[dddd]
        &{}&{}&{}
        &{\cdot}\ar@{-}[dddd]
        &{}&{}&{}
        &{\cdot}\ar@{-}[dddd]
        &{}&{}&{}
        &{\cdot}\ar@{-}[dddd]
        \\
         {}&{}&{}&{}&{}&{}&{}&{}&{}&{}
        &{}&{}&{}&{}&{}&{}&{}&{}&{}&{}
        &{}&{}&{}&{}&{}&{}&{}&{}&{}&{}
        &{}&{}&{}&{}&{}&{}&{}
        \\
        {\cdot}\ar@{-}[uurrrrr]
               \ar@{-}[dddd]
               \ar@{-}@<0.1ex>[ddrrrrr]
               \ar@{<-}[ddrrrrr]
               \ar@{-}@<-0.1ex>[ddrrrrr]
        &{}{}&{}&{} &{}&{}&{}&{}&{}&{}
        &{}&{}&{}&{}&{}&{}&{}&{}&{}
        &{\cdot}\ar@{-}[uurrrrr]
                \ar@{-}[dddd]
        &{}&{}&{}
        &{\cdot}\ar@{-}[uurrrrr]
                \ar@{-}[dddd]
        &{}&{}&{}
        &{\cdot}\ar@{-}[uurrrrr]
                \ar@{-}[dddd]
                \ar@{-}@<0.1ex>[ddrrrrr]
                \ar@{<-}[ddrrrrr]
                \ar@{-}@<-0.1ex>[ddrrrrr]
        &{}&{}&{}
        &{\cdot}\ar@{-}[uurrrrr]
                \ar@{-}[dddd]
        &{}&{}&{}&{}&{}
        \\
         {}&{}
         &{\hspace{2ex}\bullet}
         &{}&{}&{}&{}&{}&{}&{}&{}&{}&{}
         &{}
         &{}\ar@{-}@<0.1ex>[lllll]
            \ar@{->}[lllll]
            \ar@{-}@<-0.1ex>[lllll]
         &{}&{}&{}&{}&{}&{}&{}
         &{\hspace{-2ex}\bullet}
         &{}&{}&{}
         &{\hspace{-2.5ex}\bullet}                        
         &{}&{}
         &{\bullet\hspace{-1.5ex}}
         &{}&{}&{}
         &{\bullet\hspace{-1ex}}
          \ar@{.}[lllllllllllllllllllllllllllllll]
         &{}&{}&{}
        \\
        {}\ar@{-}@<0.1ex>[uurrrrr]
           \ar@{->}[uurrrrr]
           \ar@{-}@<-0.1ex>[uurrrrr]
        &{}&{}&{}&{}
        &{\cdot}
        &{}&{}&{}&{}&{}&{}&{}&{}&{}&{}
        &{}&{}&{}
        &{}\ar@{-}@<0.1ex>[uurrrrr]
           \ar@{->}[uurrrrr]
           \ar@{-}@<-0.1ex>[uurrrrr]
        &{}&{}&{}&{}
        &{\cdot}
        &{}&{}&{}
        &{\cdot}
        &{}&{}
        &{}\ar@{-}@<-0.2ex>[uurrrrr]
           \ar@{<-}@<-0.1ex>[uurrrrr]
           \ar@{-}@<0ex>[uurrrrr]
        &{\cdot}
        &{}&{}&{}
        &{\cdot}
        \\
        {}
        &{} \ar@{-}@<0.1ex>[uuuurrr]
            \ar@{->}[uuuurrr]
            \ar@{-}@<-0.1ex>[uuuurrr]
        &{}&{}&{}&{}&{}&{}&{}&{}&{}&{}
        &{}&{}&{}&{}&{}&{}&{}&{}&{}&{}
        &{}&{}
        &{} \ar@{-}@<0.1ex>[uuuurrr]
            \ar@{->}[uuuurrr]
            \ar@{-}@<-0.1ex>[uuuurrr]
        &{}&{}&{}&{}&{}&{}&{}&{}&{}&{}
        &{}&{}
        \\
        {\cdot}\ar@{-}[uurrrrr]
        &{}&{}&{}&{}&{}&{}&{}&{}&{}&{}
        &{}&{}&{}&{}&{}&{}&{}&{}
        &{\cdot}\ar@{-}[uurrrrr]
        &{}&{}&{}
        &{\cdot}\ar@{-}[uurrrrr]
        &{}&{}&{}
        &{\cdot}\ar@{-}[uurrrrr]
        &{}&{}&{}
        &{\cdot}\ar@{-}[uurrrrr]
        &{}&{}&{}&{}&{}
        \\
        {}&{}&{}&{}&{}&{}&{}&{}&{}&{}
        &{}&{}&{}&{}&{}&{}&{}&{}&{}&{}
        &{}&{}&{}&{}&{}&{}&{}&{}&{}&{}
        &{}&{}&{}&{}&{}&{}&{}
        \\
        {}&{}&{}&{}&{}&{}&{}&{}&{}&{}
        &{}&{}&{}&{}&{}&{}&{}&{}&{}&{}
        &{}&{}&{}&{}&{}&{}&{}
        &{}\ar@{-}@<0.5ex>[dd]
           \ar@{->}@<0.4ex>[dd]
           \ar@{-}@<0.3ex>[dd]
        &{}&{}&{}&{}&{}&{}&{}&{}&{}
        \\
        {}&{}&{}&{}&{}&{}&{}&{}&{}&{}
        &{}&{}&{}&{}&{}&{}&{}&{}&{}&{}
        &{}&{}&{}&{}&{}&{}&{}&{}&{}&{}
        &{}&{}&{}&{}&{}&{}&{}
        \\
        {}&{}&{}&{}&{}&{}&{}&{}&{}&{}
        &{}&{}&{}&{}&{}&{}&{}&{}&{}&{}
        &{}&{}&{}&{}&{}&{}&{}&{}&{}&{}
        &{}&{}&{}&{}&{}&{}&{}
        \\
        {}&{}&{}&{}&{}&{}&{}&{}&{}&{}
        &{}&{}&{}&{}&{}&{}&{}&{}&{}&{}
        &{}&{}&{}&{}&{}&{}&{}&{}&{}&{}
        &{}&{}&{}&{}&{}&{}&{}
        \\
         {}&{}&{}&{}&{}&{}&{}&{}&{}&{}
        &{}&{}&{}&{}&{}&{}&{}&{}&{}&{}
        &{}
        &{\bullet}\ar@{-}[rrrrrrrrrrrr]
                  \ar@{-}@<-0.1ex>[rrrrrrrrrrrr]
                  \ar@{.}@<-0.7ex>[uuuuuuuuu]
        &{}&{}&{}
        &{\bullet}\ar@{.}@<-0.7ex>[uuuuuuuuu]
        &{}&{}&{}
        &{\bullet}\ar@{.}@<-0.7ex>[uuuuuuuuu]
        &{}&{}&{}
        &{\bullet}\ar@{.}@<-0.7ex>[uuuuuuuuu]
        &{}&{}&{}
       }
$$
\caption{\emph{La recta proyectiva y la familia universal}}
\label{moduli}                           
\end{figure}
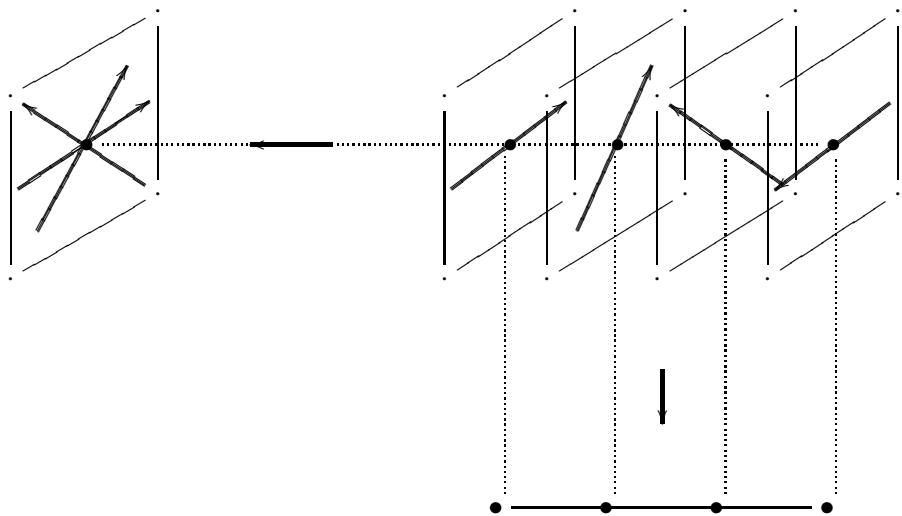

Ilustramos con el cl\'asico ejemplo de los espacios proyectivos, m\'as particularmente la recta proyectiva $\mathbb{P}^1$, que es el conjunto de rectas por el origen en el plano. Este conjunto puede considerarse como el intervalo real $[0, \, \pi]$ con los extremos identificados, a cada $x \in [0, \, \pi]$ le corresponde la recta de pendiente $x$. Se construye un conjunto $C$ en el espacio, $C \subset \mathbb{P}^1 \times \mathbb{R}^2$, munido de la proyecci\'on $C \longrightarrow \mathbb{P}^1$, poniendo sobre cada $x$ la recta de pendiente $x$ (Figura \ref{moduli}). As\'{\i}, cuando $x$ avanza a partir de $x = 0$ hacia la derecha, las rectas van girando (``modulan'') hasta llegar a $x = \pi$, que es la misma recta de partida. La continuidad de este movimiento est\'a dada matem\'aticamente por la topolog\'{\i}a del intervalo $[0, \, \pi]$ con los extremos identificados, y la del espacio $[0, \, \pi] \times \mathbb{R}^2$.

Consideremos la categor\'{\i}a $Top$ de los espacios topol\'ogicos y las aplicaciones continuas y el funtor $F:Top^{op}\to Set$ tal que $F(X)$ es el conjunto de familias continuas indexadas por $X$ de rectas del plano pasando por el origen, es decir,  
$$
F(X) = \{p: V \longrightarrow X,\; V\subseteq X\times\mathbb{R}^2, \; V_x = p^{-1}(x) = \text{una recta por el origen} \}
$$ 

$F$ act\'ua en las flechas por medio del producto fibrado: para cada aplicaci\'on continua $f:Z\to X$, $F(f)(V)=f^*(V)$, 
$$
\xymatrix@R=3ex@C=4ex
         {
           f^*(V) \ar[r]\ar[d] & V \ar[d]
           & f^*(V)_z = V_{f(z)}
           \\
           Z \ar[r]^f & X  
         }
$$

Este funtor es representable por el par 
$\mathbb{P}^1 \in Top$, 
$(C \longrightarrow \mathbb{P}^1)  \in  F(\mathbb{P}^1)$, que es la \emph{familia universal}. Esto significa que para cada familia $V\to X \in F(X)$, existe una \'unica funci\'on continua 
$f:X\to\mathbb{P}^1$ tal que $V = f^*(C)$. Se tiene un producto fibrado:
$$
\xymatrix@R=3ex@C=6ex
         {
           V \ar[r]\ar[d]_{\forall \;} & C \ar[d]
           & V_x = C_{f(x)}
           \\
           X \ar@{-->}[r]^{\exists ! f} & \mathbb{P}^1  
         }
$$
La recta $V_x$ se identifica con la recta de pendiente $f(x)$, y las rectas de la familia $V$ indexada por $X$ var\'{\i}an (``modulan'') continuamente seg\'un la funci\'on $f$.

\vspace{1ex}

Por el principio de dualidad, el lema de Yoneda es tambi\'en v\'alido para funtores covariantes $hom(C,-) : \mathcal{X}  \to  Set$, $F: \mathcal{X} \to Set$, que definen propiedades universales del otro lado de la dualidad. En realidad, es el mismo lema pero enunciado en la categor\'{\i}a dual. En el caso contravariante el morfismo \'unico ``llega'' a la soluci\'on, en el caso covariante el morfismo \'unico ``sale'' de la soluci\'on.

\section*{La teor\'{i}a de categor\'{\i}as abstractas}

{\bf [1955-1957]} {\bf Categor\'{i}as abelianas}. Alexander Grothendieck desarrolla el \'algebra homol\'ogica del Cartan-Eilenberg en categor\'{i}as abstractas definidas por propiedades categ\'oricas, llamadas \emph{categor\'{i}as abelianas} con las propiedades AB3, AB4, AB5, AB6. Esta publicaci\'on  \cite[1957]{G}\footnote{El trabajo ya estaba escrito en 1955, carta a Serre del 4 de junio 1955, y le\'{\i}do en un encuentro del grupo Bourbaki, respuesta de Serre del 13 de julio 1955, v\'ease \mbox{\emph{Correspondance Grothendieck-Serre},} Soci\'et\'e Math\'ematique de France, 2001.} engloba adem\'as de las categor\'{\i}as de m\'odulos a las categor\'{\i}as de haces de m\'{o}dulos. Se emplea por primera vez un sistem\'atico uso de propiedades universales en las definiciones y sus duales.

\vspace{1ex}

{ \bf[1958]} {\bf Funtores adjuntos}. Daniel M. Kan, al utilizar la teor\'{\i}a de categor\'{i}as en un trabajo sobre la homotop\'{i}a \cite[1958]{K}, introduce y desarrolla en abstracto el concepto de \emph{funtores adjuntos}, teniendo en mente la realizaci\'on geom\'etrica de un conjunto simplicial. 

Un par de funtores
$ 
\xymatrix
         {
           \mathcal{X} \ar@<1ex>[r]^G 
         & \mathcal{Y} \ar@<1ex>[l]^F 
         }
$ se dicen \emph{adjuntos}, denotado $F \dashv G$, si se tiene una correspondencia biun\'{\i}voca entre flechas expresada en el siguiente diagrama:
\begin{equation} \label{adjuntos}
\hspace{2ex} 
\forall \, X \in \mathcal{X}, \; Y \in \mathcal{Y},
\hspace{5ex} 
\frac{F(Y) \to X}{Y \to G(X)}
\end{equation}  

La teor\'{\i}a de categor\'{\i}as abstractas comienza a despegarse de un desarrollo abstracto del \'algebra homol\'ogica, y se  desarrolla en categor\'{\i}as que ya no son abelianas, inspirada en otras \'areas de la matem\'atica.

\vspace{1ex} 
   
{\bf [1958]} {\bf Categor\'{i}as internas}. Charles Ehresmann (1905-1979) introduce las \mbox{\emph{categor\'{i}as internas}} a partir de categr\'{i}as diferenciales y topol\'ogicas \mbox{\cite[1958]{EH}.} Llega al concepto de categor\'{i}a como generalizaci\'on del concepto de \emph{grupoide}, en particular, de grupoide diferenciable, es decir, cuyos conjuntos de objetos y de flechas son una variedad diferenciable. Por otro lado, as\'{\i} como un grupo es un grupoide con un solo objeto, quitando las inversas, concibe una categor\'{\i}a como la estructura tal que se reduce a un monoide cuando tiene un solo objeto, \emph{un grupo es a un grupoide como un monoide es a una categor\'{i}a.} 

Ehresmann desarrolla la teor\'{\i}a de categor\'{\i}as abstractas y sus nociones y resultados b\'asicos fundamentales, pero, contrariamente a la notaci\'on adoptada por Eilenberg y Mac Lane, lo hace con un \emph{lenguage desprovisto de objetos}. Como consecuencia, sus trabajos resultan ilegibles y han sido ignorados y frecuentemente redescubiertos por otros autores. 

\vspace{1ex}

{\bf[1960]} {\bf SGA1: Funtores pro-representables. Categor\'{i}as fibradas.} Grothendieck introduce en \cite[1960/61]{SGA1} estos conceptos fundamentales, y desarrolla en particular una \emph{teor\'{i}a del descenso}. Fundamentos que van mucho m\'as all\'a de sus aplicaciones a la geometr\'{\i}a algebraica. Un concepto clave en la noci\'on de categor\'{\i}a fibrada es el de \emph{morfismo cartesiano}, definido por una propiedad universal. Un \emph{clivaje} consiste en una elecci\'on de morfismos cartesianos, pero Grothendieck omite el clivaje como parte de la estructura de categor\'{\i}a fibrada. Al respecto, hace el siguiente comentario: \emph{Es por otro lado probable, al contrario del uso todav\'{i}a hoy preponderante debido a viejas maneras de pensar, que terminar\'{a} por adoptarse en los problemas universales no poner el acento en \emph{una} soluci\'{o}n elegida de una vez por todas, sino poner todas las soluciones en un pie de igualdad}\footnote{SGA1 Expos\'{e} VI: ``Il est d'ailleurs probable que, contrairement a l'usage encore pr\'{e}pond\'{e}rant maintenent, li\'{e} \`{a} d'anciennes habitudes de pens\'{e}e, il finira par s'av\'{e}rer plus commode dans les probl\'{e}mes universels, de ne pas mettre l'accent sur \emph{une} solution suppos\'{e}e choisie une fois pour toutes, mais de mettre toutes les solutions sur un pied d'\'egalit\'{e}''. \cite[p\'{a}g. 194]{SGA1}}.   

\vspace{1ex}

{\bf[1964]} {\bf El primer libro}. En 1964 aparece un hermoso librito \cite[1964]{F} del cual aprendimos mucho sobre categor\'{\i}as los estudiantes como yo, m\'as acostumbrados a los libros que a los art\'{\i}culos originales.

\vspace{1ex} 
 
{\bf[1967]} {\bf Categor\'{i}as de fracciones}. El desarrollo por P. Gabriel y M. Zisman \cite[1967]{GZ} de las categor\'{i}as de fracciones est\'a, como el concepto de funtor adjunto, inspirado en la teor\'{i}a de homotop\'{i}a.

\section*{Los topos}

{\bf [1962-1964]} {\bf Topos de Grothendieck y ``cambio de base''.} La monumental obra SGA4 de Grothendieck (con M. Artin y J-L. Verdier) \mbox{\cite[1963/64]{SGA4},} precedida por \mbox{\cite[1962]{A},} inicia la teor\'{i}a de topos al estudiar las categor\'{i}as de haces de conjuntos. Los topos aparecen como ``generalizaci\'on'' de espacio topol\'ogico. Se considera esencial la noci\'on de morfismo $\mathcal{E}\to\mathcal{F}$ de topos, que corresponde a la noci\'on de funci\'on continua, y que incorpora expl\'{\i}citamente la imagen inversa. Los topos de Grothendieck $\mathcal{E}$ son definidos por un sitio peque\~{n}o $\mathcal{C}$ (noci\'on que generaliza el reticulado de abiertos de un espacio topol\'ogico) en la categor\'{\i}a $Set$ de los conjuntos, \mbox{$\mathcal{E} = Sh(\mathcal{C})$,} y est\'an munidos de un morfismo can\'onico (y \'unico) $\mathcal{E} \to Set$. El verdadero objeto de estudio es la 2-categor\'{i}a $\mathcal{T}op$ de topos sobre $Set$, y los \emph{cambios de base}, que son sus morfismos.

\vspace{1ex}

{\bf [1962-1970]} {\bf Topos elementales.} Los topos elementales surgen como generalizaci\'on de la categor\'{i}a de los conjuntos. En el periodo 1962-1970, en su programa de fundamentar la matem\'atica con el concepto de categor\'{\i}a, \mbox{F. William Lawvere} publica una serie de art\'{\i}culos donde descubre e investiga propiedades de la categor\'{\i}a de conjuntos expresables por enunciados elementales del lenguage categ\'orico, y que son v\'alidos tambi\'en en las categor\'{\i}as de haces. Estos desarrollos culminan con la noci\'on de \emph{topos elemental} en un trabajo conjunto con Myles Tierney, presentado en el Congreso Internacional de Matem\'aticos de 1970 \cite[1970]{L}. Se considera al topos como un universo $\mathcal{S}$, aislado. Se muestra que es posible trabajar ``internamente'' en un topos $\mathcal{S}$. Lawvere introduce conectivas proposicionales por \emph{propiedad universal} y cuantificadores como \emph{funtores adjuntos} asociados a una proyecci\'on $\pi: X\times Y\to Y$: 
$$
\xymatrix
         {
          \exists_x\dashv \pi_x^{-1}\dashv\forall_x:
         }
\xymatrix@C=9ex    
         {                   
           {Partes}(Y) \ar[r]^{\pi_x^{-1}}                           
         & {Partes}(X \times Y)
                       \ar@<-1.5ex>@/_2ex/[l]_{\forall_x}
                        \ar@<1.5ex>@/^2ex/[l]_{\exists_x}
         }   
$$

Considerando $Partes$ como extensiones de f\'ormulas con las correspondientes variables libres, las correspondencias dadas por la adjunci\'on, ver diagrama (\ref{adjuntos}), resultan ser las siguientes: 
$$\frac{\phi(y)\Rightarrow\forall_x\varphi(x,y)}{\pi_x^{-1}\phi(y)\Rightarrow \varphi(x,y)}, \hspace{5ex}\frac{\exists_x\varphi(x,y)\Rightarrow\phi(y)}{\varphi(x,y)\Rightarrow\pi_x^{-1}\phi(y)},$$
donde ``$\Rightarrow$'' denota la relaci\'on de contenido entre subobjetos.
Se obtienen exactamente la reglas de la \emph{l\'ogica intuicionista}, que no son otra cosa entonces que una instancia del concepto matem\'atico de funtores adjuntos. 

El trabajo de Lawvere tiene un sustrato filos\'ofico inspirado en el materialismo que considera la conciencia como parte de la materia. As\'{i}, para Lawvere, la l\'ogica est\'a contenida en la matem\'atica, invirtiendo la relaci\'on establecida por Bertrand Russell.  

\vspace{1ex}

{\bf [1970-1972]} 
En varias conferencias impartidas en 1970-72, nunca publicadas pero ampliamente difundidas, Andr\'e Joyal hace aportes fundamentales a la teor\'{\i}a de topos. En particular:

\vspace{1ex}

{\bf Sem\'antica de Kripke-Joyal.} 
La idea es aplicar la sem\'antica de Tarski para la interpretaci\'on de las f\'ormulas l\'ogicas en la categor\'{\i}a $Set$ de los conjuntos, pero esta vez con la interpretaci\'on en un topos de Grothendieck arbritrario $\mathcal{E}$. Joyal logra describir expl\'{\i}citamente la sem\'antica resultante, se tiene una ``Sheaf semantics'' que puede ser utilizada para validar enunciados internos en topos particulares. Esta sem\'antica tiene como caso particular la sem\'antica de Kripke, y frecuentemente es referida como ``Kripke-Joyal semantics''. 

\vspace{1ex}

{\bf Topos relativos a un topos de base.} Joyal hace una s\'{\i}ntesis de Lawvere y Grothendieck. Utilizando la l\'ogica del topos, Joyal considera sitios internos en un topos $\mathcal{S}$, y muestra c\'omo construir un topos de haces (internos) 
$Sh(\mathcal{C})$, munido de un morfismo can\'onico al topos (``de base'') 
$\mathcal{S}$. Son los $\mathcal{S}$-topos. Joyal demuestra que todo morfismo de topos de \mbox{Grothendieck} $\mathcal{E} \to \mathcal{S}$ es de esa forma. Como $\mathcal{S}$ es arbritrario, dado un morfismo cualquiera entre topos de Grothendieck 
$\mathcal{E}\to\mathcal{F}$, $\mathcal{E}$ puede ser considerado como un $\mathcal{F}$-topos. Esta es la culminaci\'on del concepto de ``cambio de base''.

\section*{Categor\'{\i}as enriquecidas, Bicategor\'{\i}as}

{\bf [1966-1970]} {\bf Categor\'{i}as enriquecidas}. Primero Eilenberg y Max Kelly (1930-2007) en \cite[1966]{EK} y luego Eduardo J. Dubuc \cite[1970]{D} establecen y desarrollan las nociones de \emph{categor\'{i}as cerradas}, \emph{categor\'{i}as monoidales cerradas}, y la teor\'{i}a de \emph{categor\'{i}as enriquecidas}. Toda categor\'{\i}a $\mathcal{X}$ esta munida del funtor 
\mbox{$hom: \mathcal{X}^{op} \times \mathcal{X} \longrightarrow Set$} que a cada par de objetos le asigna el conjunto $hom(X, Y)$ de flechas entre $X$ e $Y$. Una categor\'{\i}a monoidal cerrada $\mathcal{V}$ es una categor\'{\i}a que tiene la estructura suficiente para ser receptora de un funtor $hom$, y una categor\'{\i}a $\mathcal{X}$ enriquecida es una categor\'{\i}a munida de un funtor $hom: \mathcal{X}^{op} \times \mathcal{X} \longrightarrow \mathcal{V}$. Por ejemplo, si $\mathcal{V}$ es la categor\'{\i}a de los grupos abelianos, las categor\'{\i}as enriquecidas son las \emph{categor\'{\i}as aditivas}, cuyos conjuntos $hom(X, Y)$ son un grupo commutativo. 

\vspace{1ex}

{\bf [1969]} {\bf Bicategor\'{\i}as.} 
Si $\mathcal{V}$ es la categor\'{\i}a de categor\'{\i}as, las categor\'{\i}as enriquecidas son las \emph{2-categor\'{\i}as}, los conjuntos $hom(X, Y)$ son una categor\'{\i}a cuyos objetos son las flechas $X \to Y$, y cuyas flechas, llamadas \emph{2-celdas}, son flechas entre flechas. El primer ejemplo es la categor\'{\i}a de categor\'{\i}as, cuyas flechas son los funtores, y cuyas 2-celdas son las transformaciones naturales, ver diagrama (\ref{Cat}). Los conjuntos pueden considerarse como categor\'{\i}as cuyas \'unicas flechas son las identidades, son las \emph{0-categor\'{\i}as}, en ellas solo tengo la relaci\'on de igualdad $X = Y$ entre objetos. En las categor\'{\i}as usuales o \emph{1-categor\'{\i}as} puedo comparar objetos y tengo la relaci\'on de \emph{isomorfismo} $X \cong Y$, pero solo tengo la relaci\'on de igualdad $f = g$ entre flechas. En las \emph{2-categor\'{\i}as} puedo comparar objetos de dos formas, con la relaci\'on de \emph{equivalencia} $X \sim Y$, y con la relaci\'on de isomorfismo $X \cong Y$, y puedo comparar flechas con la relaci\'on de isomorfismo $f \cong g$,  pero solo tengo la relaci\'on de igualdad $\alpha = \beta$ entre 2-celdas. Est\'a claro c\'omo esto puede seguir hasta el infinito. As\'{\i}, una \emph{3-categor\'{\i}a} es una categor\'{\i}a cuyos $hom$ son 2-categor\'{\i}as, y en ella se pueden comparar 2-celdas por medio de \emph{3-celdas}, etc. etc.  

\vspace{1ex}

 0-cat: Conjunto: $X=Y$.

 1-cat: Categor\'{i}a: $X\to Y$,\; $X\cong Y$,\; $f=g$.

 2-Categor\'{\i}a:  
$
X
\xymatrix@C=9ex@R=2.4ex
    {
     {} \ar@<1.6ex>[r]^{f}
        \ar@{}@<-1.4ex>[r]^{\!\!\alpha\,\Downarrow}
        \ar@<-1.3ex>[r]_{g}
   & {}
    }
Y\,,
$     
$X\sim Y$,\; $X\cong Y$,\; $f\cong g$,\; $\alpha=\beta$.

$\cdots $

\vspace{1ex}

Consideremos una estructura, como las 2-categor\'{\i}as, que tiene objetos, flechas y cuyos conjuntos $hom(X,Y)$ son categor\'{\i}as, es decir, tiene 2-celdas. Entonces, el diagrama original (\ref{diagrama}) con el que comenzamos esta nota, al poder comparar flechas, puede expresarse:
$$
\xymatrix
        {
         h \cong gf
        }
\hspace{20ex}
\xymatrix@R=12pt@C=8pt
         {
           X  \ar[rr]^f  \ar[rd]_h   
         & {} \ar @{}[d]|(.35){\cong} 
         &  Y \ar[ld]^g
         \\
         & Z  
         } 
$$

As\'{\i}, en una tal estructura, puedo relajar las ecuaciones de identidad y la asociatividad y poner como parte de la estructura isomorfismos
\mbox{$f id_X \cong f \cong id_Y f$, $h(gf) \cong (hg)f$.} As\'{\i} se tiene una  
estructura (que no es una categor\'{\i}a), introducida y desarrollada por B\'{e}nabou \cite[1969]{B} bajo el nombre de \emph{bicategor\'{\i}a}. A los isomorfismos de las identidades y de la asociatividad se les imponen unas pocas ecuaciones que resultan en este caso  ``dictadas'' por el propio formalismo, que son los axiomas llamados \emph{de coherencia}. 

Hoy en d\'{\i}a las bicategor\'{\i}as se llaman \emph{2-categor\'{\i}as d\'ebiles} dado que repitiendo esta relajaci\'on de la asociatividad, al poder comparar 2-celdas se puede pedir que las identidades y la asociatividad de la composicion de 2-celdas valgan ``salvo'' isomorfismos dados por 3-celdas. Estas son las \emph{3-categor\'{\i}as d\'ebiles}. En este caso los axiomas de coherencia llenan varias p\'aginas, y en el caso de las \emph{4-categor\'{\i}as debiles} nadie los ha podido escribir. Estas especulaciones no son vac\'{\i}as de contenido, ya que est\'an \'{\i}ntimamente ligadas a la noci\'on de homotop\'{\i}a entre funciones continuas, a las homotop\'{\i}as entre homotop\'{\i}as, y as\'{\i} siguiendo. Es decir, se tiene un modelo concreto, pero no se sabe cu\'ales son los axiomas que cumple.

Pero hasta el a\~{n}o 1970 solo se hab\'{\i}a llegado a $n = 2$, es decir, a las bicategor\'{\i}as.

\end{document}